\documentclass[11pt]{amsart}

\usepackage{amssymb,latexsym}

\usepackage{enumerate}

\makeatletter

\@namedef{subjclassname@2010}{%

  \textup{2010} Mathematics Subject Classification}

\makeatother

\newtheorem{theorem}{Theorem}

\newtheorem{corollary}[theorem]{Corollary}

\newtheorem{lemma}[theorem]{Lemma}

\newtheorem{proposition}[theorem]{Proposition}

\theoremstyle{definition}

\newtheorem{definition}[theorem]{Definition}

\numberwithin{equation}{section}

\newcommand{\setd}[2]{\,\left\{#1\ \colon\ #2\right\}}

\newcommand{\RR}{\mathbb{R}}

\newcommand{\NN}{\mathbb{N}}

\newcommand{\E}{{\ell_2(G)\otimes C(X)}}

\begin{document}

\title[Amenable actions]{Hilbert $C^*$-modules and amenable actions}

\author[R.G. Douglas]{Ronald G. Douglas}

\address{Department of Mathematics, Texas A\&M University, College Station, TX 77840}

\email{rdouglas@math.tamu.edu}

\author[P.W. Nowak]{Piotr W. Nowak}

\address{Department of Mathematics, Texas A\&M University, College Station, TX 77840}

\email{pnowak@math.tamu.edu}
\date{}

\begin{abstract}
We study actions of discrete groups on Hilbert $C^*$-modules induced from topological 
actions on compact Hausdorff spaces. We show non-amenability of actions of non-amenable
and non-a-T-menable groups, provided there exists a quasi-invariant probability measure which is sufficiently
close to being invariant.
\end{abstract}



\maketitle

The notion of topological amenability  of group actions has found many applications in recent years,
particularly in index theory. Yu proved \cite{yu} that the coarse Baum-Connes conjecture 
and the Novikov conjecture hold for groups which 
satisfy property A, a weak version of amenability. 
Property A turned out to be equivalent to existence of a topologically amenable 
action on some compact space \cite{higson-roe} and to exactness of the
reduced group $C^*$-algebra $C^*_r(G)$ \cite{guentner-kaminker,ozawa}. 
Because of the interest
of finding counterexamples to the the above conjectures it is natural
to study conditions which would imply non-amenability of topological actions.

Given a topological action of a non-amenable group on a compact space, the existence of
a finite invariant probability measure implies that the action is not topologically amenable 
(see Definition \ref{definition  : amenable action} and the remarks following it).
However, apart from this fact there are practically no results which would give sufficient 
conditions for
non-amenability of an action unless one assumes the existence of an invariant probability measure 
for the action.
In this paper we study the situation in which we are given a topological 
action on a compact space $X$ and a probability measure $\nu$ such that
the action of $G$ preserves the measure class. This means the translate
$g^*\nu$ and $\nu$ are absolutely continuous with respect to each other and,
in particular, the Radon-Nikodym derivatives $dg^*\nu/d\nu$ are defined almost everywhere
for every $g\in G$. The general idea is that if there is a probability measure for the action which
is sufficiently close to being invariant, then we can still prove non-amenability of 
actions using this probability measure.

Our first result is that if the
Radon-Nikodym derivatives of the translated measures satisfy 
some global integrability conditions,  then a
topologically amenable action gives rise to a proper, affine isometric action on a Hilbert space.
The latter property, known as a-T-menability or the Haagerup property, was defined
by Gromov \cite{gromov}.   As a consequence we get our first result,
namely that for groups which do not admit such actions, e.g. groups
with property (T), our condition on the Radon-Nikodym derivatives implies that
the action of $G$ cannot be topologically amenable. 

Our second result
is that if a non-amenable group $G$ acts via measure class preserving homeomorphisms
and the probability measure satisfies a certain metric condition then the action is not amenable. The condition
is expressed in terms of an inequality between the bottom of the positive spectrum 
of the discrete Laplacian on $G$ and the average Hellinger distance between $\nu$
and its translates by generators. The Hellinger distance is a bounded metric on the space of probability measures which quantifies how far the probability measure is from being invariant.
In the latter case the distance between $\nu$ and its translates is always zero. 
In the last section we discuss some examples and applications.

The main tool that we use is the fact that the action of a group $G$
on a compact space $X$ gives a linear representation of $G$ into
the   group of non-adjointable, norm-bounded, linear isometries
of the Hilbert $C^*$-module $\ell_2(G)\otimes C(X)$.
If additionally we equip $X$ with a quasi-invariant probability measure $\nu$, then we can use 
the larger module $\ell_2(G)\otimes L_{\infty}(X,\nu)$. The fact that $G$ preserves the
class of the probability measure $\nu$ allows us to overcome the non-adjointability of 
the above isometric representation and to apply Hilbert space techniques
to analyze the action and related representations.

\subsection*{Acknowledgements} We would like to thank the referee for carefully reading the manuscript and 
suggesting many improvements.

\section{ Hilbert $C^*$-modules and unitary representations}

Let $G$ be a group generated by a finite, symmetric set $S$ (i.e., $S=S^{-1}$) and let $\vert \,\cdot\,\vert:G\to\RR$
denote the associated word length function.  The word length metric on $G$ is the left-invariant
metric $d(g,h)=\vert g^{-1}h\vert$.
Let $X$ be a compact, Hausdorff space equipped with an action of
$G$ by homeomorphisms, $g\mapsto \Phi_g$. We denote the induced action of $G$ on 
$f\in C(X)$ by
automorphisms by
$$g*f(x)=f\left(\Phi_{g^{-1}}(x)\right),$$
where $f\in C(X)$.
$G$ also has a natural action on itself by left translations which induces 
the left regular representation denoted
$$g\cdot \xi_h=\xi_{g^{-1}h}$$
for $\xi\in \ell_2(G)$.

\subsection{The $G$-regular representation on Hilbert $C^*$-modules} \label{subsection : G-Hilbert modules}
For a group $G$ we will consider the following linear representations on a Hilbert  $C^*$-module.
Let the action of $G$ on a compact topological space $X$ be given.
Consider the linear space
$$\mathbb{F}=\setd{\xi:G\to C(X)}{\xi_g=0 \text{ for all but finitely many } g}.$$
Equip $\mathbb{F}$ with the inner product $\langle\, \cdot\,,\cdot\,\rangle_{C(X)}:\mathbb{F}\to C(X)$ given by
taking the regular scalar product defined for $x\in X$,
$$\langle \xi,\eta\rangle_{C(X)}(x)=\sum_{g\in G}\xi_g(x)\eta_g(x).$$
Finally, complete the resulting space in the norm 
$\Vert v\Vert_{\E}=\Vert \langle v,v\rangle_{C(X)}\Vert^{1/2}$.
The resulting space is a Hilbert $C^*$-module $\ell_2(G)\otimes C(X)$.
An analogous construction can also be done after replacing $C(X)$ with $L_{\infty}(X,\nu)$ 
for a probability measure $\nu$.
A standard reference on this material is Lance's book \cite{lance}. 

Given a Hilbert $C^*$-module $\mathcal{E}$, denote by $\operatorname{Iso}(\mathcal{E})$
the group of linear isomorphisms which preserve the norm but which are not necessarily 
adjointable.
Define a  representation $L:G\to \operatorname{Iso}(\E)$ by 
setting
$$(L_g\xi)_h(x)=\xi_{g^{-1}h}\left(\Phi_{g^{-1}}(x)\right),$$
for all $\xi\in\mathbb{F}$, $g,h\in G$ and $x\in X$ and extend to linear operators on 
$\E$.
We abbreviate $L_g\xi=g*g\cdot \xi$. The action $L_g$ is the diagonal action on $\ell_2(G)\otimes C(X)$. 
(Note that the order of applying $g$ does not matter,
since the actions $\cdot$ and $*$ commute). 
This representation satisfies 
\begin{equation}\label{equation : G-unitary representation}
\langle L_g \xi,L_g\eta\rangle_{C(X)}= g* \langle\, g\cdot \xi, g\cdot \eta\rangle_{C(X)}\end{equation}
for all $\xi,\eta\in \E$ and $g\in G$. 
Note that the operators $L_g$
are linear and bounded in norm but are not adjointable operators on the Hilbert module. However,
we can use them to construct a unitary representation of $G$.

\subsection{Unitary representations induced by $G$-unitary representations}
Consider the Hilbert $C^*$-module, $\ell_2(G)\otimes L_{\infty}(X,\nu)$. The module
$\ell_2(G)\otimes C(X)$ is a submodule of $\ell_2(G)\otimes L_{\infty}(X,\nu)$
and the above representation extends.
We introduce the following scalar product on ${\ell_2(G)\otimes L_{\infty}(X,\nu)}$:
\begin{equation}\label{equation : scalar product}
\langle \xi,\eta \rangle= \int_{{X}} \langle \xi,\eta\rangle_{L_{\infty}(X,\nu)}(x)\ {d\nu},
\end{equation}
where $\langle\xi,\eta\rangle_{L_{\infty}(X,\nu)}\in L_{\infty}(X,\nu)$ is defined analogously to $\langle \xi,\eta\rangle_{C(X)}\in C({X})$.
This turns the space ${\ell_2(G)\otimes L_{\infty}(X,\nu)}$ into a pre-Hilbert space and we obtain a Hilbert space $\mathcal{H}$
by completion. Denote
$$\rho_g=\dfrac{dg^*\nu}{d\nu},$$
so that 
$$\int_X g*f\ \rho_g d\nu=\int_X f\ d\nu$$
for every measurable function $f:X\to \RR$.
A standing assumption in this paper is that the the Radon-Nikodym derivatives
$\rho_s$, $s\in S$,  exist and are elements of $L_{\infty}(X,\nu)$. Define
$$\pi_g={\rho_g}^{1/2}\ L_g.$$
Since the Radon-Nikodym derivatives
are  elements of $L_{\infty}(X,\nu)$, each $\pi_g$ is a bounded linear operator on the
Hilbert module ${\ell_2(G)\otimes L_{\infty}(X,\nu)}$. 
Moreover, the Radon-Nikodym derivatives satisfy the cocycle condition
\begin{equation}\label{equation : cocycle conditon for Radon-Nikodym derivatives}
{\rho_{gh}}={\rho_{g}}\,g*{\rho_{h}},
\end{equation}
which guarantees that $\pi_{gh}=\pi_g\pi_h$.
Moreover, one can easily check that
$$\langle \pi_g\xi,\pi_g\eta \rangle= \langle \xi,\eta\rangle.$$
Thus  the extension of each $\pi_g$ to $\mathcal{H}$ (also denoted $\pi_g$) is a unitary
operator and we obtain a unitary representation $\pi$ of $G$ on $\mathcal{H}$. 

Unitary representations as above are often used in the context of measurable cocycles, 
see for
instance \cite{zimmer-ergodic}.

\subsection{Topologically amenable actions}
Topological amenability of homeomorphic actions on compact spaces was defined in \cite{delaroche-renault}
and was modeled on Zimmer's definition of measurable amenable actions \cite{zimmer}.
\begin{definition}\label{definition  : amenable action}
Let $X$ be a compact topological space on which $G$ acts by homemorphisms. The action 
is topologically amenable
if for every $\varepsilon>0$ there exists $\xi\in \mathbb{F}$ such that 
\begin{enumerate}
\renewcommand{\labelenumi}{(\alph{enumi})}
\item $\xi_g\ge 0$ for every $g\in G$,
\item$\langle \xi,\xi\rangle_{C(X)}=1_X$,
\item $\sup_{x\in X} \left(1-\frac{1}{\#S}\sum_{s\in S}\langle \xi, L_s\xi\rangle_{C(X)}(x)\right)\le \varepsilon$.
\end{enumerate}
\end{definition}
Amenability of an action of $G$ does not depend 
on the choice of the (finite) set of generators since we can express the new generators as finite
products of the old generators.
If $X$ is a single point, then the definition reduces to that of amenability of $G$.
Another way to phrase amenability of an action is to say that the groupoid of the action of 
$G$ on $X$
is amenable, see \cite{delaroche-renault}. This condition can also be rephrased in terms
of isoperimetric inequalities with coefficients in a $G$-$C^*$-algebra \cite{nowak-isoactions}.

It was proved in \cite{kuhn} that for measurably amenable ergodic actions the representation
$\pi$ is weakly contained in the regular representation.

It is well-known that if there exists a $G$-invariant mean, that is, a continuous linear positive,
and $G$-invariant functional on $C(X)$ (which, by the Riesz representation theorem, corresponds to an invariant probability measure on $X$),
then the action of $G$ on $X$ is amenable if and only if $G$ is amenable. Given such a mean
an appropriate averaging procedure applied to the functions $\xi:G\to C(X)$ as in Definition \ref{definition  : amenable action}
gives corresponding functions $\widetilde{\xi}:G\to \RR$ which satisfy $\Vert \widetilde{\xi}-s\cdot \widetilde{\xi}\Vert\to 0$. 

\section{ Amenable actions and a-T-menability}
In this section we will give conditions on the Radon-Nikodym derivatives
which will imply the non-amenability of actions by groups not having the Haagerup property.
\subsection{Affine isometric actions}
An affine isometric action of a group $G$ on a Banach space $E$ is given by 
$$A_gv=\pi_gv+b_g,$$
where $\pi:G\to \operatorname{Iso}(E)$ is a  representation in the linear isometry group
of $E$ and $b:G\to E$ satisfies the 
cocycle condition 
$$b_{gh}=\pi_gb_h+b_g,$$
for $g,h\in G$.
The action is called metrically proper if for every $v\in E$ we have $\lim_{\vert g\vert\to\infty} \Vert A_gv\Vert =\infty$, which is equivalent to 
$$\lim_{\vert g\vert\to\infty}\Vert b_g\Vert =\infty.$$

\begin{definition}\cite{gromov}
A group which admits a metrically proper affine isometric action on a Hilbert space
is said to be \emph{a-T-menable} or to have the \emph{Haagerup property}
\end{definition}
See \cite{cherix-et-al} for a detailed account of a-T-menability.
Let 
$$\overline{\rho}(x)=\sup_{g\in G}\rho_g(x),$$
$$\underline{\rho}(x)=\inf_{g\in G}\rho_g(x),$$
Both $\underline{\rho}$ and $\overline{\rho}$ are $\nu$-measurable since
$G$ is countable
and we have $\underline{\rho}(x)\le 1\le \overline{\rho}(x)$ for all $x\in X$.
Thus $\underline{\rho}$ is automatically an element of $L_{\infty}(X,\nu)$.

\begin{theorem}\label{theorem : firs main}
Let $G$ be a finitely generated group. Assume that $G$
acts by homeomorphisms on a compact Hausdorff space $X$ and that there is 
a probability measure $\nu$ on $X$ such that at least one of the following conditions holds
\begin{enumerate}
\item $\overline{\rho}\in L_1(X,\nu)$, \label{enumerate : supremum}  
\item $\int_X\ \underline{\rho}(x)\ d\nu>0$.\label{enumerate : infimum}
\end{enumerate}
If the action is topologically amenable then the group admits a proper affine isometric action on 
a Hilbert space.
\end{theorem}
If the probability measure $\nu$ is invariant, then $\rho_g= 1_X$ for every $g\in G$,
and the above conditions are trivially satisfied.
In particular, we get
\begin{corollary}
With the conditions of the above theorem, if $G$ is not a-T-menable
then the action is not amenable. 
\end{corollary}\label{corollary : main a-T-menable}
Examples of groups which are not a-T-menable include groups which have property (T)
or relative property (T). See the monographs \cite{cherix-et-al} and
\cite{bekka-delaharpe-valette}.

Boundedness conditions related  to but stronger than (\ref{enumerate : supremum}) 
were studied by Greenleaf 
\cite{greenleaf}, Feldman and Moore \cite{feldman-moore}. They
showed that if the Radon-Nikodym derivatives are globally bounded,
which in our case translates to $\overline{\rho}\in L_{\infty}(X,\nu)$,
 then there
exists an equivalent invariant probability measure. Later, Zimmer \cite{zimmer-cohomology}
showed that if 
$\sup_{g\in G}\rho_g(x)<\infty$  for almost every $x\in X$ and, in addition the action is 
ergodic, then again there is an equivalent invariant probability measure. For general actions it is not known
whether the conditions in Theorem \ref{theorem : firs main} imply the existence of an 
equivalent invariant probability measure. Condition (\ref{enumerate : infimum}), however, has not 
been studied in this context.
Note that by the cocycle property we have that
$$\underline{\rho}(x)>0\ \ \ \Longleftrightarrow\ \ \  \sup_{g\in G} g^{-1}*\rho_g(x)<\infty.$$

The next statement shows that the existence of an amenable action on a space
gives an affine isometric action on a Hilbert module, which is in addition assumed to be
metrically proper in a certain stronger sense.

\begin{proposition}\label{proposition : affine action on a Hilbert module}
Let $G$ be a finitely generated group acting amenably on a compact space $X$.
Then $G$ admits an affine isometric action on a Hilbert module  $\mathcal{E}$ over $C(X)$
with a cocycle $b:G\to \mathcal{E}$ such that the following functional inequality holds 
$$\phi(\vert g\vert) 1_X\ \le\ \langle b_g,b_g\rangle_{C(X)}\ \le\ K \vert g\vert^2 1_X$$
for some $K>0$ and some nondecreasing $\phi:[0,\infty)\to[0,\infty)$, $\lim_{t\to\infty} \phi(t)=\infty$.
\end{proposition}

\begin{proof}
We use the cocycle construction as in \cite{bekka-cherix-valette}.
Let $\xi_n$ be as in Definition \ref{definition  : amenable action} with $\varepsilon=1/n$. 
Let $R_n\in \RR$ be such that $\langle \xi_n,L_g\xi_n\rangle_{C(X)}=0$ whenever 
$\vert g\vert \ge R_n$. Take the Hilbert module 
$$\mathcal{E}=\bigoplus_{n\in \NN}\ell_2(G)\otimes C(X)$$
and define the representation $T:G\to \operatorname{Iso}(\mathcal{E})$ by 
$$T_g=\bigoplus_n L_g.$$
Then define a function $b:G\to \mathcal{E}$ by
$$b_g=\bigoplus_{n\in\NN} L_g\xi_n-\xi_n.$$
One can easily check that $b_g$ is a cocycle for $T$ and we will now estimate 
its norm.
For any generator $s\in S$, 
\begin{eqnarray*}
\langle b_s,b_s\rangle_{C(X)}&=&\sum_{n\in \NN} \langle L_s\xi_n-\xi_n,L_s\xi_n-\xi_n\rangle_{C(X)}\\
&\le &\sum_{n\in\NN} \dfrac{1}{n^2} 1_X.
\end{eqnarray*}
Thus letting $K=\sum_{n\in \NN} \frac{1}{n^2}$ we obtain the upper bound
for $g\in S$. The bound for general $g\in G$ follows by writing $g$ as a word in 
generators and applying the upper bound for each of these.

Since each $\xi_n$ is finitely supported we have that 
$$\langle L_g\xi_n-\xi_n,L_g\xi_n-\xi_n\rangle_{C(X)}=2\cdot 1_X$$ 
whenever $\vert g\vert\ge R_n$. For $g\in G$ let $\phi(\vert g\vert)$ 
be the largest $n$
for which $\vert g\vert\ge R_n$. Then we have
\begin{eqnarray*}
\langle b_g,b_g\rangle_{C(X)}&=&\sum_{n\in \NN} \langle L_g\xi_n-\xi_n,L_g\xi_n-\xi_n\rangle_{C(X)}\\
&\ge &2\sum_{n=1}^{\phi(\vert g\vert)} 1_X\\
&=&2\phi(\vert g\vert) 1_X. 
\end{eqnarray*}
It is not hard to see that
$\lim_{t\to\infty}\phi(t)=\infty$ and we thus obtain the lower bound $2\phi$.
\end{proof}

A crucial property of $\overline{\rho}$ and $\underline{\rho}$  is the
following invariance.
\begin{lemma}\label{lemma : identities for rho's}
The following identities hold for any $g\in G$:
\begin{enumerate}
\item ${\rho_g}^{1/2}\left(g*\overline{\rho}^{\,1/2}\right)=\overline{\rho}^{\,1/2}$,
\item ${\rho_g}^{1/2}\left(g*\underline{\rho}^{\,1/2}\right)=\underline{\rho}^{\,1/2}.$
\end{enumerate}
\end{lemma}
\begin{proof}
We will prove (1); (2) is completely analogous.
For any fixed $g\in G$ we have
\begin{eqnarray*}
{\rho_g}(x)^{1/2}g*\overline{\rho}(x)^{1/2}&=&\sup_{h\in G}{\rho_g}(x)^{1/2}\ g*{\rho_h}(x)^{1/2}\\
&=&\sup_{h\in G} { \rho_{gh}}(x)^{1/2}\\
&=&\overline{\rho}(x)^{1/2},
\end{eqnarray*}
for all $x\in X$.
\end{proof}

\begin{proof}[Proof of Theorem \ref{theorem : firs main}]
First we prove the assertion under the assumption (\ref{enumerate : infimum}).
Let $b:G\to \bigoplus \ell_2(G)\otimes C(X)$ be a cocycle for the action 
as in Proposition \ref{proposition : affine action on a Hilbert module}.
It can also be viewed as a cocycle $b:G\to\bigoplus\ell_2(G)\otimes L_{\infty}(X,\nu)$
for the same representation, viewed now as a representation on 
$\bigoplus\ell_2(G)\otimes L_{\infty}(X,\nu)$.
Define a function $\underline{b}:G\to \bigoplus\ell_2(X)\otimes L_{\infty}(X,\nu)$ 
by
$$\underline{b}_g=\underline{\rho}^{1/2}\,b_g.$$
Then $\underline{b}$ is a cocycle for the unitary representation $U_g=\bigoplus\pi_g$ on the
Hilbert space $\mathcal{H}$ defined via the scalar product $\langle \oplus\xi_n,\oplus \eta_n\rangle=\sum \langle \xi_n,\eta_n\rangle$, where 
the summands are scalar products defined by equation (\ref{equation : scalar product}). Indeed, we have
$$U_g=\bigoplus {\rho_g}^{1/2} L_g={\rho_g}^{1/2} T_g$$
and, by Lemma \ref{lemma : identities for rho's},
\begin{eqnarray*}
U_g\underline{b}_h+\underline{b}_g&=&{\rho_g}^{1/2}T_g\underline{\rho}^{1/2}b_h+\underline{\rho}^{1/2}b_g\\
&=&{\rho_g}^{1/2} (g*\underline{\rho}^{1/2})T_gb_h+\underline{\rho}^{1/2}b_g\\
&=&\underline{\rho}^{1/2}\left(T_gb_h+b_g\right)\\
&=&\underline{b}_{gh}.
\end{eqnarray*}
We now have
\begin{eqnarray*}
\Vert \underline{b}_g\Vert^2&=&\int_X\langle \underline{b}_g,\underline{b}_g\rangle(x)\ d\nu\\
&=&\int_X\underline{\rho}(x)\langle b_g,b_g\rangle(x)\ d\nu.
\end{eqnarray*}
Applying the functional inequalities from Proposition \ref{proposition : affine action on a Hilbert module}  we obtain
$$\phi(\vert g\vert)\left( \int_X\underline{\rho}(x)\ d\nu \right)\le  \Vert \underline{b}_g\Vert^2\le K\vert g\vert^2 \left(\int_X\underline{\rho}(x)\ d\nu\right).$$
Thus if $\int_X\underline{\rho}(x)\ d\nu =\underline{C}>0$ then the affine isometric
action 
$$\underline{A}_gv=U_gv+\underline{b}_g$$
on $\mathcal{H}$ is well-defined and metrically proper since $\underline{C}\phi(t)\to\infty$ as $t\to\infty$.

To prove the assertion assuming (\ref{enumerate : supremum}) we 
take $\overline{b}=\overline{\rho}^{1/2}b$. Similarly, $\overline{b}$ is a cocycle for $\pi$.
In this case we need an additional argument.
This is because given a vector $v\in\ell_2(G)\otimes C(X)$, the vector $\overline{\rho}^{1/2}b_g$
will not be an element of the Hilbert module $\ell_2(G)\otimes L_{\infty}(X,\nu)$ unless $\overline{\rho}$ is bounded (which is exactly 
what we are trying to avoid). However,
if $\overline{\rho}\in L_1(X,\nu)$ then 
$\overline{\rho}^{1/2}b$ is an element of a Hilbert space 
$\bigoplus\left(\ell_2(G)\otimes L_2(X,\nu)\right)$.
Repeating the above argument for the cocycle $\overline{b}$ we obtain 
$$\phi(\vert g\vert)\left( \int_X\overline{\rho}(x)\ d\nu\right)\le \Vert \overline{b}_g
\Vert^2\le
K\vert g\vert \left(\int_X\overline{\rho}(x)\ d\nu\right).$$
If $\int_X\overline{\rho}(x)\ d\nu <\infty$ then the isometric
affine action $$\overline{A}_gv=U_gv+\overline{b}_g$$
is well-defined and metrically proper.
\end{proof}

\section{ Actions of non-amenable groups}
In this section we will give a different condition for the non-amenability of an action of
a non-amenable group. As mentioned earlier, if there is an invariant probability measure for an 
action of such a group, it follows easily that the action is not topologically amenable. However 
we are
interested in the situation in which the  probability measure is only quasi-invariant.
Similar ideas were used in \cite{nevo}, however, with different motivations.

\subsection{The Hellinger distance for probability measures}
Given probability measures $\mu_1$ and $\mu_2$, both absolutely continuous
with respect to the probability measure $\nu$ on $X$, we consider the formula
$$H(\mu_1,\mu_2)=\left(\dfrac{1}{2}\int\left(\sqrt{\dfrac{d\mu_1}{d\nu}}-\sqrt{\dfrac{d\mu_2}{d\nu}}\right)^2\ d\nu \right)^{1/2}.$$
$H$ does not depend on the choice of the dominating measure $\nu$ and is
known as the Hellinger distance between probability distributions
\cite{pollard,vandervaart}.
We can also write 
$$H(\mu_1,\mu_2)=\left(1-A(\mu_1,\mu_2) \right)^{1/2},$$
where the quantity 
$$A(\mu_1,\mu_2)=\int_X\sqrt{\dfrac{d\mu_1}{d\nu}\,\dfrac{d\mu_2}{d\nu}}\ d\nu$$
is referred to as the Hellinger affinity. The Fubini theorem applied to $A$ gives one of the fundamental properties of the Hellinger metric, namely its behavior with respect to product
measures. The Hellinger metric satisfies the following inequalities \cite[page 61]{pollard} with respect to the $L_1$-metric:
$$H(\mu_1,\mu_2)^2\le \Vert \mu_1-\mu_2\Vert_{L_1}\le H(\mu_1,\mu_2).$$
The Hellinger distance is used in asymptotic statistics and in quantum mechanics
(see \cite{pollard,vandervaart}).
Note that  $H(\mu_1,\mu_2)=0$ if and only if $\mu_1=\mu_2$ 
and $H(\mu_1,\mu_2)=1$ if and only if $\mu_1$ and $\mu_2$  are singular. 

\subsection{Spectrum of the Laplacian}
Let $G$ be an infinite group generated by a finite set $S$. The bottom of the spectrum of the discrete Laplace operator
on the Cayley graph $X=(V,E)=\mathcal{G}(G,S)$ is defined via the variational expression 
$$\lambda_1=\inf_{f\in \ell_2(G)}\dfrac{\langle df,df\rangle}{\langle f,f\rangle}
=\dfrac{\sum_{s\in S,\,g\in G} \vert f_g-f_{s^{-1}g}\vert^2}{\sum_{g\in G}\vert f_g\vert^2},$$
where $d:\ell_2(V)\to \ell_2(E)$ is defined by $df(x,y)=f(y)-f(x)$ for an edge $(x,y)\in E$.
The group $G$ is amenable if $\lambda_1=0$ for any Cayley graph of $G$ and if $\lambda_1>0$, then it gives 
a sort of measure of how non-amenable $G$ is. The constant $\lambda_1$ is
closely related to the Cheeger constant of $G$,  defined by
$$h=\inf\setd{\dfrac{\#\partial F}{\#F}}{F\subset G \text{ is finite }}.$$
In particular, $h>0$ if and only if $\lambda_1>0$.
See e.g. \cite{chavel,lubotzky} for background.

\begin{theorem}\label{theorem : hellinger distance}
Let $G$ be a non-amenable group generated by a finite set $S$ and
$X$ be a compact Hausdorff space equipped with a probability measure $\nu$.
Let $G$ act on $X$ by homeomorphisms which preserve the measure class of $\nu$.
If
\begin{equation}\label{equation : Hellinger dist condition}
\dfrac{1}{\#S}\sum_{s\in S}H(\nu,s^*\nu)^2 <\dfrac{\lambda_1}{2},\end{equation}
then the action of $G$ on ${X}$ is not amenable. 
\end{theorem}

\begin{proof}
Let 
$$\beta=\dfrac{1}{\#S}\sum_{s\in S}\int_X\rho_s(x)^{1/2}\ d\nu.$$
Then $\frac{1}{\# S}\sum_{s\in S}H^2(\nu,s^*\nu)=1-\beta$.
Assume the action of $G$ on ${X}$ is amenable and consider $\xi$ for the given $\varepsilon>0$
as in Definition \ref{definition  : amenable action}. 
Then $\xi$ 
satisfies
$$\langle \xi,\xi\rangle=\int_{X}\langle \xi,\xi\rangle_{C(X)}(x)\ {d\nu}=\int_{X} 1\ {d\nu}=1.$$
Moreover, since the Radon-Nikodym derivatives are
positive, we have
\begin{eqnarray*}
\dfrac{1}{\#S}\sum_{s\in S}\langle \xi, \pi_s\xi\rangle&=&\dfrac{1}{\#S}\sum_{s\in S} \int_{X} {\rho_s}(x)^{1/2} \langle \xi, L_s\xi\rangle_{C(X)}(x)\ {d\nu}\\
&\ge&\dfrac{1}{\#S}\sum_{s\in S}\int_{{X}}{\rho_s}(x)^{1/2} (1-\varepsilon)\ {d\nu}\\
&\ge& \beta(1-\varepsilon).
\end{eqnarray*}
Since each such $\xi$ is finitely supported, for each $\xi$ there is an $R>0$ such that $\langle L_g\xi,\xi\rangle_{C(X)}=0$ for all $g\in G$ satisfying
$\vert g\vert\ge R$. Consequently for such $\vert g\vert \ge R$ we have
$$\langle \pi_g\xi,\xi\rangle= \int_{{X}}{\rho_g}^{1/2}(x) \langle \xi_x,L_g\xi_x\rangle_{C(X)} (x) {d\nu}=0.$$
If we set
$$\psi_g=\langle \pi_g\xi,\xi\rangle,$$
then $\psi$ is a finitely supported positive definite function on $G$.
Since $\psi$ is finitely supported, it defines a bounded convolution operator $T$ 
on $\ell_2(G)$ and thus an element of  
$C^*_r(G)$, the reduced $C^*$-algebra of $G$.
By the positive definiteness of $\psi$, the operator $T$ is positive and there exists a square root
$Q\in C^*_r(G)$ such that 
$$Q^*Q=T.$$
We now define  $\eta\in \ell_2(G)$ by setting
$$\eta_g=(Q1_e)g,$$
where $1_e$ is the point mass at  $e$.
Thus $\eta$ is the column labeled by $e$ in the matrix representation of $Q$.
By the definition of $\eta$ we have
$$\langle \eta,g\cdot \eta\rangle=\psi_g.$$
Here $\langle\,\cdot\,,\,\cdot\,\rangle$ denotes the standard scalar product in $\ell_2(G)$. 
We conclude that $\Vert \eta\Vert=1$ 
and 
\begin{eqnarray*}
\dfrac{1}{\#S}\sum_{s\in S}\langle \eta,s\cdot\eta\rangle&=&\dfrac{1}{\#S}\sum_{s\in S}\psi(s)\\
&=&\dfrac{1}{\#S}\sum_{s\in S}\langle \pi_s\xi,\xi\rangle\\
&\ge& \beta(1-\varepsilon).
\end{eqnarray*}
Together with the definition of the isoperimetric constant $\lambda_1$, we obtain
the inequality
\begin{eqnarray*}
\lambda_1&\le &\langle d\eta,d\eta\rangle\\
&=&\dfrac{2}{\#S}\sum_{s\in S}\left(1-\langle \eta,s\cdot\eta\rangle\right)\\
&\le&2\left( 1-\beta(1-\varepsilon)\right)
\end{eqnarray*}
for every $\varepsilon>0$ and consequently, ${\lambda_1}\le 2(1-\beta)$, which is a contradiction.
\end{proof}

\section{Concluding remarks}

\subsection{Actions which are close to isometric actions}

We can study the question of amenability of actions in a setting
where we require an action to be close to 
an isometric action. See for example \cite{fisher-margulis} for such a study in the context of 
rigidity.
Given an action of a group $G$ on a compact manifold, if 
there exists a subset $U$ of positive volume on which the
action of $G$ distorts the volume by a uniformly small amount, 
then Theorem \ref{theorem : firs main} applies.
To keep calculations simple we consider the case 
of the circle.

Denote by $S^1$ the circle and by $\operatorname{Diff}^1_+(S^1)$ 
the group of orientation preserving $C^1$-diffeomorphisms of $S^1$. Consider
the ``distance at $x$" given by restricting the $r$-uniform distance to $x\in S^1$:
$$d_x(\varphi,\phi)=d_{S^1}(\varphi(x),\phi(x))+\vert D\varphi(x)-D\phi(x)\vert,$$
where $\varphi,\phi\in \operatorname{Diff}^1_+(S^1)$ and $D$ denotes the derivative.
Assume that $G$ acts on $S^1$ by diffeomorphisms, with the action given by a homomorphism 
$\varphi:G\to\operatorname{Diff}^1_+(S^1)$.
Assume also that $G$ is not a-T-menable. 
If there exists  $U\subseteq S^1$ and an isometric action 
$\phi:G\to \operatorname{Diff}^1_+(S^1)$ 
such that $\sup_{x\in U} d_x(\varphi_g,\phi_g)\le C< 1$ for all 
$g\in G$, then the action $\varphi_g$ 
is not topologically amenable.

Indeed, in that case, $\vert 1-D\varphi_g(x)\vert=\vert D\phi_g(x)-D\varphi_g(x)\vert\le C<1$,
since $\phi$ is an isometry.
This implies $\inf_{g\in G} D\varphi_g(x)> 0$. Since $D\varphi_g=\rho_g$ for every $x\in U$, the claim
follows from Theorem \ref{theorem : firs main}.
The above discussion generalizes easily to piecewise smooth homeomorphism.

\subsection{Non-amenable actions of the free group on $S^1$}

A similar fact  holds when we consider a non-amenable group. In that case we can 
restrict our attention to the generators but we have to compare the distances on the whole
circle. Let $G$ be a non-amenable, finitely generated group acting on $S^1$ by $C^1$ 
diffeomorphisms. If the generators of $G$ are sufficiently close to isometries in the sense
that
$$1- \dfrac{1}{\#S}\sum_{s\in S}\int_X\sqrt{\vert D\varphi_s\vert}\ d\nu  
<\dfrac{\lambda_1}{2},$$
then the action is not topologically amenable. 

An explicit example can be constructed as follows.
Introduce the $C^1$-topology on $\operatorname{Diff}^1_+(S^1)$ by the metric
$$d(\varphi,\phi)=\sup_{x\in S^1}d_{S^1}(\varphi(x),\phi(x))+\sup_{x\in S^1}\vert D\varphi(x)-D\phi(x)\vert.$$
The $C^1$ topology turns $\operatorname{Diff}^1_+(S^1)$ into a Baire space
\cite{palis-demelo}.
Adapting the transversality argument from \cite[Proposition 4.5]{ghys} we see that 
any generic (in the sense
of Baire's category theorem) pair
of diffeomorphisms in $\operatorname{Diff}^1_+(S^1)$ generates a free group.
Consider any pair of isometries 
$(i_1,i_2)\in \operatorname{Diff}^1_+(S^1)\times  \operatorname{Diff}^1_+(S^1)$.
Arbitrarily close
to the pair $(i_1,i_2)$ there exists a pair $(q_1,q_2)$ which generates a free group.  
Clearly, $d(q_k,i_k)\le\varepsilon$ implies
$$ \sup_{x\in S^1}\vert Dq_k(x)-1\vert \le \varepsilon.$$
However, in this case the $Dq_k$, $k=1,2$, are the Radon-Nikodym derivatives and thus we conclude 
that they satisfy the conditions of Theorem \ref{theorem : hellinger distance} if $\varepsilon$ is
sufficiently small.  

\subsection*{Acknowledgements}

The first author was supported by NSF grant DMS-0600865. The second author was supported by NSF grant DMS-0900874.


\begin{thebibliography}{22}

\bibitem{delaroche-renault}
C.~Anantharaman-Delaroche, J.~Renault, \emph{
Amenable groupoids.} With a foreword by Georges Skandalis and Appendix B by E. Germain. 
Monographies de L'Enseignement MathŽmatique 36. L'Enseignement MathŽmatique, Geneva, 
2000.

\bibitem{bekka-cherix-valette}
{M.E.B.~Bekka, P.-A.~Cherix,  A.~Valette},  \emph{Proper affine isometric actions of 
amenable groups.} Novikov conjectures, index theorems and rigidity, Vol. 2 (Oberwolfach, 1993), 
1--4, London Math. Soc. Lecture Note Ser., 227, Cambridge Univ. Press, Cambridge, 1995.

\bibitem{bekka-delaharpe-valette}
{B.~Bekka, P.~de la Harpe, A.~Valette}, \emph{Kazhdan's property ($T$).} 
New Mathematical Monographs, 11. Cambridge University Press, Cambridge, 2008.

\bibitem{chavel}
{I.~Chavel}, \emph{Topics in isoperimetric inequalities.} Geometry, spectral theory, groups, and dynamics, 89--110, Contemp. Math., 387, Amer. Math. Soc., Providence, RI, 2005.

\bibitem{cherix-et-al}
{P.-A.~Cherix, M.~Cowling, P.~Jolissaint, P.~Julg, A.~Valette}, \emph{Groups 
with the Haagerup property. Gromov's a-T-menability.} Progress in Mathematics, 197. Birkh\"{a}user Verlag, Basel, 2001. 

\bibitem{feldman-moore}
{J.~Feldman, C.C~Moore}, \emph{Ergodic equivalence relations, cohomology, and von Neumann algebras. II.}  Trans. Amer. Math. Soc.  234  (1977), no. 2, 325--359.

\bibitem{fisher-margulis}
{D.~Fisher, G.~Margulis}, \emph{Almost isometric actions, property (T), and local rigidity.} Invent. Math. 162 (2005), no. 1, 19--80.

\bibitem{ghys}
{E.~Ghys}, \emph{Groups acting on the circle.} Enseign. Math. (2) 47 (2001), no. 3-4, 329--407.

\bibitem{greenleaf}
{F.P.~Greenleaf}
\emph{Invariant means on topological groups and their applications.}
Van Nostrand Mathematical Studies, No. 16 Van Nostrand Reinhold Co., New York-Toronto, Ont.-London 1969.

\bibitem{gromov}
{M.~Gromov}, \emph{Asymptotic invariants of infinite groups.} Geometric group theory, Vol. 2 (Sussex, 1991), 1--295, London Math. Soc. Lecture Note Ser., 182, Cambridge Univ. Press, Cambridge, 1993.

\bibitem{guentner-kaminker}
{E.~Guentner, J.~Kaminker}, \emph{Exactness and the Novikov conjecture.} 
Topology 41 (2002), no. 2, 411--418. 

\bibitem{higson-kasparov}
{N.~Higson, G.~Kasparov}, \emph{$E$-theory and $KK$-theory for groups which act properly and isometrically on Hilbert space.} Invent. Math. 144 (2001), no. 1, 23--74.


\bibitem{higson-roe}
{N.~Higson, J.~Roe}, \emph{Amenable group actions and the Novikov conjecture.}
J. Reine Angew. Math. 519 (2000), 143--153. 

\bibitem{kuhn}
{M.G.~K\"{u}hn}, \emph{Amenable actions and weak containment of certain representations of discrete groups.}  Proc. Amer. Math. Soc.  122  (1994),  no. 3, 751--757.


\bibitem{lance}
{E.C.~Lance}, \emph{Hilbert $C\sp *$-modules.
A toolkit for operator algebraists.} London Mathematical Society Lecture Note Series, 210. Cambridge University Press, Cambridge, 1995.



\bibitem{lubotzky}
{A.~Lubotzky}, \emph{Discrete groups, expanding graphs and invariant measures.} With an appendix by Jonathan D. Rogawski. Progress in Mathematics, 125. Birkh\"{a}user Verlag, Basel, 1994.

\bibitem{nevo}
{A.~Nevo}, \emph{The spectral theory of amenable actions and invariants of 
discrete groups.} Geom. Dedicata 100 (2003), 187--218.

\bibitem{nowak-isoactions}
{P.W.~Nowak}, \emph{Isoperimetry of group actions.}  Adv. Math.  219  (2008),  no. 1, 1--26.

\bibitem{ozawa}
{N.~Ozawa}, \emph{Amenable actions and exactness for discrete groups.} C. R. Acad. Sci. Paris SŽr. I Math. 330 (2000), no. 8, 691--695.

\bibitem{palis-demelo}
{J.~Palis Jr.,W.~de Melo}, \emph{Geometric theory of dynamical systems. An introduction.}Translated from the Portuguese by A. K. Manning. Springer-Verlag, New York-Berlin, 1982.

\bibitem{pollard}
{D.~Pollard}, \emph{A user's guide to measure theoretic probability.} Cambridge Series in
Statistical and Probabilistic Mathematics, 8. Cambridge University Press, Cambridge, 2002.


\bibitem{vandervaart}
{W.~van der Vaart}, \emph{Asymptotic statistics.} Cambridge Series in Statistical and Probabilistic Mathematics, 3. Cambridge University Press, Cambridge, 1998. 
\bibitem{yu}
{G.~Yu}, \emph{The coarse Baum-Connes conjecture for spaces which admit a uniform embedding into Hilbert space.}  Invent. Math.  139  (2000),  no. 1, 201--240.

\bibitem{zimmer}
{R.J.~Zimmer}, \emph{Amenable ergodic actions, hyperfinite factors, and Poincar\'{e} flows.}
Bull. Amer. Math. Soc. 83 (1977), no. 5, 1078--1080. 

\bibitem{zimmer-cohomology}
{R.J.~Zimmer}, \emph{On the cohomology of ergodic group actions.}  Israel J. Math.  35  (1980), no. 4, 289--300.

\bibitem{zimmer-ergodic}
{R.J.~Zimmer}, \emph{Ergodic theory and semisimple groups.} Monographs in Mathematics, 81. Birkh\"{a}user Verlag, Basel, 1984.

\end{thebibliography}
\end{document}